\documentclass{amsart}

\usepackage[utf8]{inputenc}
\usepackage[english]{babel}
\usepackage{amsthm,amssymb,amsmath,amsrefs,pst-all,graphicx}

\newtheorem{lem}{Lemma}
\newtheorem*{classic}{Euclidean half-space theorem}
\newtheorem*{heisenberg}{Vertical half-space theorem in Heisenberg space}

\title{On the half-space theorem\\ for minimal surfaces in Heisenberg space}
\author{Tristan Alex}

\begin{document}
\begin{abstract}
 We propose a simple proof of the vertical half-space theorem for Heisenberg space.
\end{abstract}

\maketitle 
\section{Introduction}
A half-space theorem states that the only properly immersed minimal surface which is contained
in a half-space is a parallel translate of the boundary of the half-space, namely a plane. Hoffman and Meeks
first proved it for $\mathbb R^3$ (\cite{HM}). It fails in $\mathbb R^n$ or $\mathbb H^n$, 
$n\geq4$.

In recent years, there has been increased interest in homogeneous $3$-manifolds (cf. Abresh/Rosenberg \cite{AR1}, 
Hauswirth/Rosenberg/Spruck \cite{HRS}). The original proof of Hoffman and Meeks also works 
in Heisenberg space $\operatorname{Nil}_3$ with respect to \emph{umbrellas}, which are 
the exponential image of a horizontal tangent plane (\cite{AR2}).
Daniel and Hauswirth extended the theorem to vertical half-spaces of Heisenberg space, where 
vertical planes are defined as the inverse image of a straight line in the base
of the Riemannian fibration $\operatorname{Nil}_3
\to\mathbb R^2$ (\cite{DH}).
\begin{heisenberg}[Daniel/Hauswirth 2009]
 Let $S$ be a properly immersed minimal surface in Heisenberg space. If $S$ lies to one side of a vertical
plane $P$, then $S$ is a plane parallel to $P$.
\end{heisenberg}

Essential for the proof of half-space theorems is the existence of a family of catenoids or
generalized catenoids. Their existence  
is simple to establish in spaces where they can be represented as ODE solutions. For instance,
horizontal umbrellas in Heisenberg space are invariant under rotations around the vertical axis, 
so they lead to an ODE.
 However, the lack
of rotations about horizontal axes means that the existence of analogues of a horizontal catenoid amounts to establishing
true PDE solutions. Daniel and Hauswirth use a Weierstra\ss{}-type representation to reduce this problem to a
system of ODEs. Only after solving a period problem they obtain the desired family of surfaces. 

In the present paper we introduce a simpler approach: we take a coordinate model of Heisenberg space and consider
coordinate surfaces of revolution. Provided we can choose a family of surfaces whose mean curvature normal points
into the half-space, the original maximum principle argument of Hoffman and Meeks will prove the theorem. 
Our
approach is based on an idea by Bergner (\cite{B}), who 
generalized the classical half-space theorem to surfaces with negative Gaussian curvature
such that the principal curvatures satisfy an inequality, and Earp / Toubiana (\cite{ET}), who consider special Weingarten surfaces with mean curvature satisfying an inequality.

It is an open problem to prove a vertical half-space theorem for $\operatorname{PSL}_2(\mathbb R)$, where it
would apply to surfaces whose mean curvature is the so-called \emph{magic number} $H_0=1/2$, namely the 
limiting value of the mean curvature of large spheres. Here, it would state that
surfaces with mean curvature $H_0=1/2$ lying on the mean convex side of a horocylinder can only be
horocylinders, that
is, the inverse image of a horocycle of the fibration $\operatorname{PSL}_2(\mathbb R)\to\mathbb H^2$. 
Our strategy could also work there. However, so far
we have not been successful to establish the desired family of generalized catenoids with $H\leq H_0$.

I would like to thank my advisor Karsten Gro\ss{}e-Brauckmann for his help.

\section{The Euclidean half-space theorem}
\begin{classic}[Hoffman/Meeks 1990]
 A properly immersed minimal surface $S$ in $\mathbb R^3$ lying in a half-space $H$ is a plane
parallel to $P=\partial H$.
\end{classic}
\begin{proof}
 By the standard maximum principle we can assume $\operatorname{dist}(S,P)=0$ but $S\cap P=\emptyset$.

Let $\mathcal C_r\subset\mathbb R^3\setminus H$ 
be a half catenoid with necksize $r$ and $\partial\mathcal C_r\subset P$.
By the properness of $S$, we can translate $S$ by $\varepsilon>0$ towards $\mathcal C_1$ 
such that $S$ intersects $P$ but
stays disjoint to $\partial \mathcal C_r$ for all $r\in(0,1]$.

As $r$ tends to $0$, the family of catenoids 
$\mathcal C_r$ converges to $P$ minus a point. We claim that the set $I$ of parameters for which
$\mathcal C_r$ does not intersect $S$ is open. 
Consider a catenoid $\mathcal C_{r_0}$ that does not intersect $S$.
For each $r\in(0,1)$ there exists a compact set $K$ such that the distance between $\mathcal C_r$ and $P$ is larger
than $2\varepsilon$ in the complement of $K$. We may choose $K$ in a way that this property holds
for all $r$ in a small neighbourhood of $r_0$. This implies that the distance between $S$ and all these
$\mathcal C_r$ is larger than $\varepsilon$ in the complement of $K$ (cf. figure \ref{fig1}).

However, within the compact set $K$, the distance between $S$ and $C_{r_0}$ is positive, so for all
$r$ in a (possibly smaller) neighbourhood of $r_0$, this distance is still positive.

\begin{figure}[h!]
\includegraphics[]{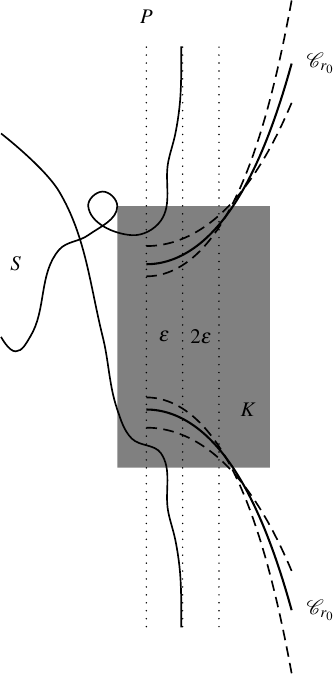}
\caption{Proof of the Euclidean half-space theorem}
\label{fig1}
\end{figure}

We conclude that in small neighbourhood of $r_0$,
\[
 \operatorname{dist}(\mathcal C_r,S)\geq\min(\operatorname{dist}(\mathcal C_r\cap K,S\cap K),
\operatorname{dist}(\mathcal C_r\cap K^c,S\cap K^c))>0,
\]
thereby proving our claim.

Therefore, the set of parameters for which $\mathcal C_r$ and $S$ do intersect
is closed, so there is a first catenoid
$\mathcal C_{r_1}$ touching $S$ at a point $p$. 
Since the boundaries of all $\mathcal C_{r}$ with $r\in(0,1]$ are disjoint from $S$, 
the touching point $p$ is an interior point, contradicting the maximum principle.
\end{proof}
\section{Coordinate surfaces of revolution}
We take the following coordinates:
\begin{align*}
 \operatorname{Nil}_3&:=(\mathbb R^3,\mathrm{d} s^2), & 
\mathrm{d} s^2&=\mathrm{d} x^2+\mathrm{d} y^2+(2\tau x\mathrm{d} y-\mathrm{d} z)^2 &
\text{with }\tau&\geq0.
\end{align*}
An orthonormal frame of the tangent space is given by
\begin{align*}
 E_1&=\partial_x, & E_2&=\partial_y+2\tau x\partial_z, & E_3&=\partial_z,
\end{align*}
and the Riemannian connection in these coordinates is determined by
\begin{equation}\label{connection}
\begin{aligned}
 \nabla_{E_1}E_2&=-\nabla_{E_2}E_1=\tau E_3,& \nabla_{E_1}E_3&=\nabla_{E_3}E_1=-\tau E_2,\\
\nabla_{E_2}E_3&=\nabla_{E_3}E_2=\tau E_1, & \nabla_{E_i}E_j&=0\text{ in all other cases.}
\end{aligned}
\end{equation}
The Heisenberg space is a Riemannian fibration $\pi\colon\mathbb R^3\to\mathbb R^2$ with vanishing base curvature.
The bundle curvature of $\operatorname{Nil}_3$ is given by $\frac12g(\nabla_{E_1}E_2-\nabla_{E_2}E_1,E_3)=\tau$ and
for $\tau=0$ we recover $\mathbb R^3$.

Let us consider a curve $c(t)=(0,t,r(t))$ in Heisenberg space with a positive function $r$ and $t\geq 0$. 
By rotating around the $y$-axis, we get an immersion
\[
 f\colon[t_0,\infty)\times[0,2\pi)\to\operatorname{Nil}_3,\quad (t,\varphi)\mapsto 
\begin{pmatrix}-r(t)\sin\varphi\\t\\r(t)\cos\varphi\end{pmatrix}.
\]

In order to apply the proof of Hoffman/Meeks, we will construct Euclidean rotational surfaces around the $y$-axis.
With the Heisenberg space metric, these rotations are not isometric, because
the $4$-dimensional isometry group of $\operatorname{Nil}_3$ contains
only translations and rotations around the vertical axis. Therefore, the mean curvature of such a surface will
depend on the angle of rotation $\varphi$. We will need to find a surface with mean curvature
vector pointing to the half-space to arrive at the desired contradiction with the maximum principle.

The tangent space of $M:=f([t_0,\infty)\times[0,2\pi))$ is spanned by
\begin{align*}
 v_1&=-r'(t)\sin\varphi E_1+E_2+(2\tau r(t)\sin\varphi+r'(t)\cos\varphi)E_3,\\
 v_2&=-r(t)\cos\varphi E_1-r(t)\sin\varphi E_3,
\end{align*}
so the inner normal of $M$ is
\[
 N=\frac{1}{W}(\sin\varphi E_1+(r'(t)+2\tau r(t)\sin\varphi\cos\varphi)E_2-\cos\varphi E_3),
\]
where $W=\sqrt{1+(2\tau r(t)\sin\varphi\cos\varphi+r'(t))^2}$.

We will now compute the first and second fundamental forms of $M$. We easily get
\begin{align*}
 G_{ij}&=\mathrm{d}s^2(v_i,v_j)\\&=\begin{pmatrix}
 \sin^2\!\varphi\, r'(t)^2+\left(2 \tau  r(t) \sin \varphi+\cos \varphi r'(t)\right)^2+1 & -2 \tau  r(t)^2 \sin ^2\varphi\\
 -2 \tau  r(t)^2 \sin ^2\varphi & r(t)^2 \\
\end{pmatrix}
\end{align*}
with determinant $\det G=r(t)^2W^2$.

The most tedious part of the calculation is the second fundamental form. We have to compute 
$B_{ij}=\mathrm{d}s^2(\nabla_{v_i}v_j,N)$. 
 To start, (\ref{connection}) gives
\begin{align*}
\nabla_{v_1}E_1&=(-2\tau^2r(t)\sin\varphi-\tau r'(t)\cos\varphi)E_2-\tau E_3, \\
\nabla_{v_1}E_2&=(2\tau^2r(t)\sin\varphi+\tau r'(t)\cos\varphi)E_1-\tau r'(t)\sin\varphi E_3,\\
\nabla_{v_1}E_3&=\tau E_1+\tau r'(t)\sin\varphi E_2.
\end{align*}
We calculate
\begin{align*}
 \nabla_{v_1}v_1=&-r''(t)\sin\varphi E_1+(2\tau r'(t)\sin\varphi+r''(t)\cos\varphi) E_3\\
&-r'(t)\sin\varphi\nabla_{v_1}E_1
+\nabla_{v_1}E_2+(2\tau r(t)\sin\varphi+r'(t)\cos\varphi)\nabla_{v_1}E_3\\
=&(-r''(t)\sin\varphi+4\tau^2r(t)\sin\varphi+2\tau r'(t)\cos\varphi)E_1\\
&+(4\tau^2r(t)r'(t)\sin^2\varphi+2\tau r'(t)^2\sin\varphi\cos\varphi)E_2\\
&+(2\tau r'(t)\sin\varphi+r''(t)\cos\varphi)E_3,
\end{align*}
and obtain the first entry of $B$ as 
\begin{align*}
 B_{11}=\frac{1}{W}\bigg(
-r''(t)+4\tau^2r(t)r'(t)^2\sin^2(\varphi)+8\tau^3r(t)^2r'(t)\sin^3(\varphi)\cos(\varphi)+4\tau^2r(t)
\sin^2(\varphi)\\
+2\tau r'(t)^3\sin(\varphi)\cos(\varphi)+4\tau^2rr'^2\sin^2(\varphi)\cos^2(\varphi)
\bigg).
\end{align*}
The other three entries arise similarly from
\begin{align*}
 \nabla_{v_2}v_1=\nabla_{v_1}v_2= &-(\tau  r(t) \sin \varphi+r'(t)\cos \varphi)E_1\\
&+(\tau  r(t) \left(2\tau  r(t) \sin\varphi\cos\varphi+r'(t)\cos (2 \varphi\right))E_2\\
&+(\tau r(t) \cos\varphi-r'(t)\sin \varphi)E_3,\\
 \nabla_{v_2}v_2=&r(t) \sin \varphi E_1-2 \tau  r(t)^2\sin\varphi\cos \varphi E_2-r(t)\cos \varphi E_3.
\end{align*}
They are
\begin{align*}
 B_{12}=B_{21}&=\frac{\tau  r(t) \left(4 \tau  r(t)r'(t)
   \sin \varphi \cos ^3\varphi+\tau ^2 r(t)^2 \sin ^2(2\varphi)+\cos (2 \varphi) r'(t)^2-1\right)}{W},\\
 B_{22}&=-\frac{r(t)
   \left(\tau  r(t) \sin (2 \varphi) \left(\tau  r(t) \sin (2 \varphi)+r'(t)\right)-1\right)}{W}.
\end{align*}
We obtain the mean curvature $H$ for our coordinate surface of revolution:
\begin{lem}\label{H}
 The mean curvature $H=H(t,\varphi)$ of $f$ is given by
\begin{align*}
 H&:=\frac12\operatorname{tr}(G^{-1}B)\\
 &=\frac{G_{22}B_{11}-G_{12}B_{21}-G_{21}B_{12}+G_{11}B_{22}}{2r(t)^2W^2}\\
&=\frac{1+r'(t)^2-r(t)r''(t)+4 \tau ^2 r(t)^2 \sin ^4\varphi+2\tau r(t)r'(t)\sin\varphi\cos\varphi}
{2 r(t)W^{3}}.
\end{align*}
\end{lem}
\section{Half-space theorem in Heisenberg space}
As expected, for $\tau=0$ Lemma \ref{H} recovers the mean curvature for rotational surfaces in Euclidean space. 
For $\tau\neq0$, the
two additional terms depending on $\varphi$ in the nominator of $H$ arise because 
the horizontal rotation is not an isometry of Heisenberg space. Our goal is to exhibit a family of  
rotational surfaces satisfying $H\leq0$ with respect to the normal $N$.

Consider the rotational surface $f_c$ given in terms of
\begin{align}\label{r}
 r_c(t):=\exp\left(\frac1c\exp(ct)\right)
\end{align}
with  $c>c_0:=4\tau^2+2\tau+1$. We claim that this surface satisfies $H\leq0$ for $t>0$. 
Indeed, the following estimate
for the denominator of $H$
holds:
\begin{align*}
 2r(t)W^{3}H&\leq 1+r_c'(t)^2-r_c(t)r_c''(t)+4\tau^2r_c(t)^2+2\tau r_c(t)r_c'(t)\\
&=1+r_c(t)^2(\exp(ct)(2\tau-c)+4\tau^2)\\
&\leq1+r_c(t)^2(4\tau^2+2\tau-c)\leq 1+4\tau^2+2\tau-c\leq0.
\end{align*}

Since we consider a rotational surface with an embedded meridian, the embeddedness of 
$M_c:=f_c([t_0,\infty)\times[0,2\pi))$ is obvious. Also, the
boundary $\partial M_c=\{\exp(1/c)\cdot(\sin\varphi,0,\cos\varphi)\colon\varphi\in[0,2\pi)\}$ is explicitly known.

It is also important to note that for each $c$ and any given $\varepsilon>0$, there exists a compact
set such that the distance between $M_c$ and the plane $\{y=0\}$ is larger than $\varepsilon$ in the complement of
this compact set.

Let us summarize the result:
\begin{lem}\label{f}
 The coordinate surface of revolution whose meridian is defined by (\ref{r}) satisfies for $c>c_0$
\begin{enumerate}
 \item $H\leq0$ with respect to the normal $N$,
 \item for $c\to\infty$, the surface $M_c$ converges uniformly to a subset of $\{y=0\}$ on compact sets,
 \item $M_c$ is properly embedded,
 \item $\partial M_c=\{\exp(1/c)\cdot(\sin\varphi,0,\cos\varphi)\colon\varphi\in[0,2\pi)\}$ for all $c$. 
\end{enumerate}
\end{lem}

Using the surfaces $M_c$, our proof of the Euclidean half-space theorem literally applies to Heisenberg space.

\bibliography{literatur}{}
\bibliographystyle{plain}

\end{document}